\documentclass{elsarticle}
\usepackage{amsmath,amssymb,arydshln}
\usepackage{microtype,csquotes,url}
\newcommand{\R}{\mathbb R}
\newcommand{\N}{\mathbb N}
\newcommand{\bq}{\begin{equation}}
\newcommand{\eq}{\end{equation}}
\newtheorem{Rem}{Remark}[section]
\newcommand\gO{\mathcal{O}}
\newcommand{\te}{{\tilde{e}}}
\newcommand{\teta}{{\tilde{\eta}}}
\usepackage{xcolor,todonotes}
\newcommand{\added}[1]{{#1}}
\newcommand{\changed}[1]{{#1}}
\makeatletter
\def\ps@pprintTitle{\let\@mkboth\@gobbletwo%
 \let\@oddhead\@empty
 \let\@evenhead\@empty
 \def\@oddfoot{\scriptsize\fontsize{8}{13}\selectfont\slshape\hskip-0em
   To be published in Computer Methods in Applied Mechanics and Engineering,
   }%
 \let\@evenfoot\@oddfoot}
\makeatother
\usepackage{fancyhdr}
\begin{document}
\begin{frontmatter}
\title{On \added{Asymptotic} Global Error Estimation and Control of Finite Difference Solutions for Semilinear Parabolic Equations}
\author{Kristian Debrabant}
\ead{debrabant@imada.sdu.dk}
\address{University of Southern Denmark, Department of Mathematics and Computer Science, Campusvej 55, 5230 Odense M, Denmark}
\author{Jens Lang}
\ead{lang@mathematik.tu-darmstadt.de}
\address{Technische Universit\"{a}t Darmstadt, Fachbereich Mathematik, Dolivostr. 15,
64293 Darmstadt, Germany}
\begin{abstract}
The aim of this paper is to extend the global error estimation and control addressed in Lang and Verwer~%\cite{LaVe2007}
[SIAM J. Sci. Comput. 29, 2007] for initial value problems to finite difference solutions of semilinear parabolic
partial differential equations. The approach presented there is combined with an estimation of
the PDE spatial truncation error \added{by Richardson extrapolation} to estimate the overall error in the computed solution.
\added{Approximations of the error transport equations for spatial and temporal global errors are derived by using asymptotic estimates that neglect higher order error terms for sufficiently small step sizes in space and time.}
\added{Asymptotic} control in a discrete $L_2$-norm is achieved through tolerance proportionality and uniform or adaptive mesh refinement. Numerical examples are used to illustrate the reliability of the estimation and control strategies.
\end{abstract}

\begin{keyword}
Numerical integration for PDEs\sep method of lines\sep
finite difference method\sep \added{asymptotic} global error estimation\sep
\added{asymptotic} global error control\sep defects and local errors\sep
tolerance proportionality
%\\MSC 2000:  Primary: 65M06, 65M15, 65M20.
%\\{1998 ACM Computing Classification System:} G.1.7, G.1.8.
\end{keyword}
\end{frontmatter}\thispagestyle{fancy}
\section{Introduction} \label{intro}
We consider \changed{semilinear parabolic partial differential equations}
\begin{multline} \label{PDEsystem} %
\partial_t u(t,x) = L(t,x)u(t,x)+g(t,x,u(t,x)),\quad t\in (0,T]\,,\; x\in\Omega\subset\R^d\,,
\end{multline}
\added{in $d\in\N$ space dimensions,} \changed{where $L$ is an elliptic operator, and assume that}
an appropriate system of boundary conditions and
the initial condition%
\bq \label{PDEinitial} %
u(0,x) = u_0(x)\,,\quad x\in\overline\Omega
\eq %
are given.
The \changed{initial boundary value problem} is assumed to be well posed and to have a unique continuous
solution $u(t,x)$. % which has sufficient regularity.

The method of lines is used to solve (\ref{PDEsystem}) numerically.
We first discretize the PDE in space by means of finite differences \added{of order $q>1$}
on a (possibly non-uniform) spatial mesh $\Omega_h$ and solve the
resulting system of ODEs using existing time integrators. For simplicity, we shall
assume that this system of time-dependent ODEs can be written in the general form %
\bq \label{ODEsystem} %
\begin{array}{rll} %
U'_h(t) &\!=\!& F_h(t,U_h(t))\,,\qquad t\in(0,T]\,,\\[0.2cm]
U_h(0) &\!=\!& U_{h,0}\,,
\end{array}
\eq %
with a unique solution vector $U_h(t)$ being a grid function on
$\Omega_h$. Let %
\begin{equation}\label{Rhdef}
R_h:\,u(t,\cdot\,) \rightarrow (R_hu)(t)
\end{equation}
be the usual restriction operator defined by
$(R_hu)(t)=(u(t,x_1),\ldots,u(t,x_N))^T$, where $x_i\in\Omega_h$ and
$N$ is the number of all mesh points. Then we take as
initial condition $U_{h,0}=R_hu(0)$.

To simplify the following derivations, we assume that $F_h$ is given by
\begin{equation}
F_h(t,U_h)=L_h(t)U_h+G_h(t,U_h)
\end{equation}
with a finite difference approximation $L_h$ of $L$, and $G_h(t,R_hu)=R_hg(t,\cdot,u(t,\cdot))$.

To solve the initial value problem (\ref{ODEsystem}), we apply a
numerical integration method \added{of order $p\ge 1$} at a certain time grid %
\bq \label{timegrid} %
0 = t_0 < t_1 < \cdots < t_n < \cdots < t_{M-1} < t_M = T\,, %
\eq %
using local control of accuracy. This yields approximations
$V_h(t_n)$ to $U_h(t_n)$, which may be calculated for other values
of $t$ by using a suitable interpolation method provided by the
integrator. The global time error is then defined by %
\bq \label{globTimeError} %
e_h(t) = V_h(t) - U_h(t)\,. %
\eq %
Numerical experiments in \cite{LaVe2007} for ODE systems have shown
that classical global error estimation based on the first
variational equation is remarkably reliable. In addition, having the
property of tolerance proportionality, that is, there exists a
linear relationship between the global time error and the local
accuracy tolerance, $e_h(t)$ can be successfully controlled by a
second run with an adjusted local tolerance. Numerous techniques
to estimate global errors are described in \cite{Sk86}. A comparison
of various adaptive grid methods for partial differential equations
and implementation issues are presented in \cite{WoSaSc98,WoSaScTh05}.

In order for the method of lines to be used efficiently, it is necessary
to take also into account the spatial discretization error.
Defining the spatial discretization error by %
\bq \label{globSpaceError} %
\eta_h(t) = U_h(t) - (R_hu)(t)\,, %
\eq %
the vector of overall global errors $E_h(t)=V_h(t)-(R_hu)(t)$ may be
written as sum of the global time and spatial error, that is, %
\bq \label{globError} %
E_h(t) = e_h(t) + \eta_h(t)\,. %
\eq %
\added{We assume that $u(t,x)$ is $(q+2)$-times continuously differentiable with respect to $x$ and
$(p+1)$-times continuously differentiable with respect to $t$. Then, with maximum step sizes
$h_{max}$ in space and $\tau_{max}=\max_{i=0,\ldots,M-1}(t_{n+1}-t_n)$ in time it holds for
the global space and time error that $\|\eta_h(t_n)\|=\gO(h_{max}^q)$ and
$\|e_h(t_n)\|=\gO(\tau_{max}^p)$, $n=1,\dots,M$, respectively.}

\added{Although a posteriori error estimates and adaptive algorithms for the efficient
solution of parabolic problems are well established (see e.g.\ \cite{LangBOOK,No1996}
and references therein), the separation of global time and spatial discretization
errors is still a challenge. First experiences to estimate and balance the spatial
discretization error and the error due to time integration of the ODEs within
the method of lines have been made by Sch\"onauer, Schnepf, and Raith \cite{ScScRa1984}.
In their control strategy, the spatial mesh is initially chosen and remains fixed.
The spatial truncation error is designated to be the level to which the local time error
must be adapted. Lawson, Berzins, and Dew \cite{LaBeDe1991} proposed to additionally
control the local time error with respect to the contribution of the existing error
from the previous time steps to the global error at the end of the next time step.
The error in time is enabled to vary in relation to the spatial discretization error,
ensuring that the method of lines with a fixed spatial mesh is being used efficiently.
A successful attempt to assess and to equilibrate the individual discretization
errors with respect to a given quantity of interest has been made by Schmich and Vexler
\cite{ScVe2008}. An adjoint linear parabolic problem has to be solved backwards
in time to derive useful error bounds, which are used to enhance the resolution in
time and space to meet a user-prescribed accuracy tolerance.}

It is the purpose of this paper to present a new \added{asymptotic} error control
strategy for the global errors $E_h(t)$, \added{based on asymptotic estimates.}
We will mainly focus on reliability. So our aim is to provide error estimates
$\tilde{E}_h(t)\approx E_h(t)$ which are not only asymptotically
exact, but also work reliably for moderate tolerances, that is for
relatively coarse discretizations. \added{Approximations of the error transport equations
for spatial and temporal global errors are derived by using asymptotic estimates that
neglect higher order error terms for sufficiently small step sizes in space and time.}
The \added{approximate} global errors are measured in discrete $L_2$-norms.
A priori bounds for the global error in such norms
are well known, see e.g.\ \cite{LaTh2005,ThomasBOOK}. However, reliable a
posteriori error estimation and efficient control of the accuracy
of the solution numerically computed to an imposed tolerance level
are still challenging. We achieve \added{asymptotic} global error control by
iteratively improving the temporal and spatial discretizations
according to \added{asymptotic} estimates of $e_h(t)$ and $\eta_h(t)$. The global time
error is estimated and controlled along the way fully described in
\cite{LaVe2007}. To estimate the global spatial error, we follow an
approach proposed in \cite{Be1988} (see also \cite{LaBeDe1991}) and
use Richardson extrapolation to set up a linearised error transport
equation. \added{Both strategies have to be combined in the right manner in
order to make sure that they work reliably. Therefore, we have developed
an appropriate control rule for the global spatial error. To control
the overall global error more efficiently, we also consider
a new fully space-time adaptive approach.}

\added{Throughout the
paper we will use the terms 'approximation' and 'estimation' in the sense
of asymptotic estimates, i.e., estimates that involve the Landau symbol $\gO$.}

\added{The outline of this paper is as follows: In Section \ref{secSpaceTimeError}, we will linearize the transport equations for the global spatial and the global time error. These contain the residual time error and the spatial truncation error, which are approximated in Sections \ref{secEstRes} and \ref{secEstTrunc}. In Section \ref{secIntFormulas} we describe the discretization formulas used to approximate the solutions of the error transport equations, as well as the strategies used to adaptively adjust the time step size and the spatial mesh in dependence on the residual time error and the spatial truncation error. Now that we have approximations to the global time and global spatial error, Section \ref{secControlRule} suggests strategies to adapt the local tolerances such that in further runs first the global time error and then the global spatial error respect some global tolerances provided by the user. Finally, numerical examples and a summary are given in Sections \ref{secNumExamples} and \ref{sec:Summary}.}

\section{Spatial and time error \label{secSpaceTimeError}}
By making use of the restriction operator $R_h$, the spatial
truncation error is defined by %
\bq \label{truncerror} %
\alpha_h(t) = (R_hu)'(t) - F_h(t,(R_hu)(t))\,. %
\eq %
From (\ref{ODEsystem}) and (\ref{truncerror}), it follows that the
global spatial error $\eta_h(t)$ representing the accumulation of
the spatial discretization error is the solution of the initial
value problem %
\bq \label{odeSpatialErr} %
\begin{array}{rll} %
\eta'_h(t) &\!=\!& F_h(t,U_h(t)) - F_h(t,(R_hu)(t)) -
\alpha_h(t)\,,\qquad
t\in (0,T]\,,\\[0.2cm]
\eta_h(0) &\!=\!& 0\,.
\end{array} %
\eq %
Assuming \changed{$F_h$} to be \added{twice} continuously differentiable, the
mean value theorem for vector functions \added{applied to $\tilde{g}(\xi)=F_h(t,(R_hu)(t)+\xi\eta_h(t))$} yields %
\bq \label{odeSpatialErrLin} %
\begin{array}{rll} %
\eta'_h(t) &\!=\!& \partial_{U_h}F_h(t,U_h(t))\,\eta_h(t) -
\alpha_h(t) +
\gO(\eta_h(t)^2),\qquad t\in(0,T],\\[0.2cm]
 \eta_h(0) &\!=\!& 0\,.
\end{array} %
\eq %
With $V_h(t)$ being the continuous extension of the numerical
approximation to (\ref{ODEsystem}), the residual time error is defined by %
\bq \label{reserror} %
r_h(t) = V'_h(t) - F_h(t,V_h(t))\,. %
\eq %
Thus the global time error $e_h(t)$ fulfills the initial value
problem %
\bq \label{odeTimeErr} %
\begin{array}{rll} %
e'_h(t) &\!=\!& F_h(t,V_h(t)) - F_h(t,U_h(t)) + r_h(t)\,,\qquad
t\in (0,T]\,,\\[0.2cm]
e_h(0) &\!=\!& 0\,.
\end{array} %
\eq %
Again, the mean value theorem yields %
\bq \label{odeTimeErrLin} %
\begin{array}{rll} %
e'_h(t) &\!=\!& \partial_{U_h}F_h(t,V_h(t))\,e_h(t) + r_h(t) +
\gO(e_h(t)^2),\qquad t\in(0,T],\\[0.2cm]
e_h(0) &\!=\!& 0\,.
\end{array} %
\eq %
Apparently, by implementing proper choices of the defects
$\alpha_h(t)$ and $r_h(t)$, solving (\ref{odeSpatialErrLin}) and
(\ref{odeTimeErrLin}) will in leading order provide approximations
to the true global error. The issue of how to approximate the
spatial truncation error and the residual time error will be
discussed in \changed{Sections \ref{secEstRes} and \ref{secEstTrunc}}.

\section{\changed{Approximation} of the residual time error \label{secEstRes}}%
\changed{The numerical approximation of the global time error $e_h(t)$ as defined in (\ref{odeTimeErrLin})
requires the construction of an appropriate nearby solution $V_h(t)$ which is used in (\ref{reserror}) to define the residual time error $r_h(t)$. The usual way is to construct an
interpolatory polynomial from the numerical solutions by using Lagrange or Hermite interpolation. The latter one exploits the
fact that with approximations $V_{h,n}:=V_h(t_n)$ at certain time points also first derivatives $F_{h,n}:=F_h(t_n,V_{h,n})$ are given. In the following we present an approach proposed in
\cite{LaVe2007} to obtain the nearby solution through piecewise cubic Hermite interpolation.
It turns out that this is useful as long as $1\le p\le 3$ with $p$ being the order of the time
integration method. One step methods of order less or equal three are quite popular in the method of lines approach, since they are easy to program and the number of the arising linear systems
is still of moderate size.}

At every subinterval $[t_n,t_{n+1}]$ we form %
\bq \label{hermitePolAnsatz} %
V_h(t)=V_{h,n}+A_n(t-t_n)+B_n(t-t_n)^2+C_n(t-t_n)^3,\qquad t_n\leq
t\leq t_{n+1}, %
\eq %
and choose the coefficients such that $V'_h(t_{n})=F_{h,n}$ and
$V'_h(t_{n+1})=F_{h,n+1}$. This
gives %
\bq \label{hermitePol} %
V_h(t_n+\theta\tau_n)=v_0(\theta)V_{h,n}+v_1(\theta)V_{h,n+1}
+\tau_nw_0(\theta)F_{h,n}+\tau_nw_1(\theta)F_{h,n+1} %
\eq %
with $0\leq\theta\leq1$, and %
\bq \label{hermitePolCoeff} %
v_0(\theta)=(1-\theta)^2(1+2\theta),\;
v_1(\theta)=\theta^2(3-2\theta),\; w_0(\theta)=(1-\theta)^2\theta,\;
w_1(\theta)=\theta^2(\theta-1), %
\eq %
which imply
\bq\label{eq:Vhtnp0.5} %
V_h(t_{n+1/2}) =
\frac12 (V_{h,n}+V_{h,n+1})+\frac{\tau_{\added{n}}}8(F_{h,n}-F_{h,n+1}) %
\eq %
and %
\bq\label{eq:Vhptnp0.5} %
V'_h(t_{n+1/2}) =
\frac{3}{2\tau_{\added{n}}} (V_{h,n+1}-V_{h,n})-\frac14 (F_{h,n}+F_{h,n+1}).%
\eq %
With \eqref{eq:Vhtnp0.5} and \eqref{eq:Vhptnp0.5} we compute from (\ref{reserror}) the residual time error halfway the step interval as
\bq %
\begin{array}{rll} \label{reshalfwaydetailed}%
\displaystyle
r_h(t_{n+1/2})&\!=\!&\frac{3}{2\tau_{\added{n}}}(V_{h,n+1}-V_{h,n})-
\frac14 (F_{h,n}+F_{h,n+1})\\[0.4cm]
&&-F_h\left(t_{n+\frac12},\frac12 (V_{h,n}+V_{h,n+1})+
\frac{\tau_{\added{n}}}8(F_{h,n}-F_{h,n+1})\right).
\end{array} %
\eq%
On the other hand, \added{assuming that $F_h$ is four times continuously differentiable with respect
to the solution,} we obtain from (\ref{reserror}) by applying the Simpson rule that
\begin{multline}\label{eq:approxresidualint}
  \int_{t_n}^{t_{n+1}}r_h(t)~dt=(V_{h,n+1}-V_{h,n})-
  \frac{\tau_{\added{n}}}6(F_{h,n}+F_{h,n+1})\\-\frac23\tau_{\added{n}} F_h\left(t_{n+\frac12},V_h(t_{n+1/2})\right)+\gO(\tau_{\added{n}}^5)
\end{multline}
and consequently
\begin{equation}\label{eq:meanvaluerhapprox}
\frac1{\tau_{\added{n}}}\int_{t_n}^{t_{n+1}}r_h(t)~dt=\frac23r_h(t_{n+1/2})+\gO(\tau_{\added{n}}^4).
\end{equation}
\added{As $r_h(t_{n+\frac12})=\gO(\tau_n^{\min\{p,4\}})$, the approximation \eqref{eq:meanvaluerhapprox} is useful as long as $p\leq3$.} Then, as in \cite[Section 2.1]{LaVe2007} we consider instead of (\ref{odeTimeErrLin}) the step size frozen
version %
\bq \label{odeTimeErrFrozen} %
\begin{array}{rll} %
\te'_h(t)&\!\!=\!\!&\partial_{U_h}F_h(t_n,V_{h,n})\,\te_h(t)+
\frac23r_h(t_{n+\frac12}),\; t\in(t_n,t_{n+1}],\;n=0,\dots,M\!-\!1,\\
\te_h(0)&\!\!=\!\!&0
\end{array} %
\eq %
to approximate the global time error $e_h(t)$.
%
%\begin{Rem}
%In contrast to the classical approach in \cite{DeLa2013,LaVe2007}, the derivation of %\eqref{odeTimeErrFrozen} chosen here does not rely on the first variational equation and hence %the local error of the one-step method.
%\hfill$\Diamond$ %
%\end{Rem}
% Jens: Scheint mir zu verwirrend hier.
%
\changed{\begin{Rem}\label{rem:cubicdefectgiveslocalerror}
When defined as above by using cubic Hermite interpolation, $r_h(t_{n+1/2})$ can also be used
to retrieve in leading order the local error $\delta_{n+1}$ at time $t_{n+1}$ of any one-step method of order $1\leq p\leq 3$ through the relation
\bq \label{reshalfway} %
r_h(t_{n+1/2}) =
\frac{3}{2}\frac{\delta_{n+1}}{\tau_n} + \gO(\tau_n^{p+1})\,, %
\eq %
(see also \cite[Section 2.2]{LaVe2007} and \cite{DeLa2013}). So controlling $r_h(t_{n+1/2})$
in a local step size procedure is equivalent to the error-per-unit-step strategy (EPUS), which
gives the favourite property of tolerance proportionality \cite{ShampineBOOK} and will also be
exploited in our numerical tests.
\hfill$\Diamond$ %
\end{Rem}}%
\added{\begin{Rem}\label{rem:ChoiceOfContinuousExtension}%
Defining the continuous extension by other means than by cubic Hermite interpolation is possible. In this case, however, the approximation \eqref{eq:meanvaluerhapprox} will in general not hold, but one could use, e.g., \eqref{eq:approxresidualint}. The advantage of \eqref{eq:meanvaluerhapprox} is that $r_h(t_{n+\frac12})$ can be efficiently used to control local time stepping as described in Section \ref{secIntFormulas}.
\hfill$\Diamond$ %
\end{Rem}}

\section{\changed{Approximation} of the spatial truncation error \label{secEstTrunc}}
An efficient strategy to estimate the spatial truncation error by
Richardson extrapolation is proposed in \cite{Be1988}. We will
adopt this approach to our setting.

Suppose we are given a second semi-discretization of the PDE system
(\ref{PDEsystem}), now with doubled local mesh sizes $2h$, %
\bq \label{ODEsystem2h} %
\begin{array}{rll} %
U'_{2h}(t) &\!=\!& F_{2h}(t,U_{2h}(t))\,,\qquad t\in(0,T]\,,\\[0.2cm]
U_{2h}(0) &\!=\!& U_{2h,0}\,.
\end{array}
\eq %
In practice, one first chooses $\Omega_{2h}$ and constructs then
$\Omega_{h}$ through uniform refinement. We assume that the solution
$U_{2h}(t)$ to the discretized PDE on the coarse mesh $\Omega_{2h}$
exists and is unique. \added{For Lipschitz continuous $F_{2h}$, this condition is fulfilled.}
We define the
restriction operator $R_{2h}^h$ from the fine grid $\Omega_h$ to the
coarse grid $\Omega_{2h}$ by the identity $R_{2h}=R_{2h}^hR_h$
\added{(where $R_h$ and $R_{2h}$ are defined by \eqref{Rhdef} on $\Omega_h$ and $\Omega_{2h}$, respectively)}
and
set %
\bq %
\eta_h^c(t)=R_{2h}^h\eta_h(t),\quad U_h^c(t)=R_{2h}^hU_h(t),
\quad V_h^c(t)=R_{2h}^hV_h(t)\,. %
\eq %
From the second assumption it follows that %
\bq \label{spatialErrAssump}%
\eta_h^c(t)=2^{-q}\eta_{2h}(t)+\gO(h^{q+1}) %
\eq %
and therefore %
\bq \label{R2hu}%
R_{2h}u(t)=\frac{2^q}{2^q-1}U_h^c(t)-\frac1{2^q-1}U_{2h}(t)+\gO(h^{q+1})\,. %
\eq %
The relation $U_h^c(t)-U_{2h}(t)=\eta_h^c(t)-\eta_{2h}(t)$ together
with (\ref{spatialErrAssump}) gives %
\bq \label{UhcU2h}%
U_h^c(t)-U_{2h}(t)=\frac{1-2^q}{2^q}\eta_{2h}(t)+\gO(h^{q+1})\,. %
\eq %
The spatial truncation error on the coarse mesh $\Omega_{2h}$ is
analogously to (\ref{truncerror}) defined as %
\bq %
\alpha_{2h}(t) = (R_{2h}u)'(t)-F_{2h}(t,R_{2h}u(t))\,. %
\eq %
Substituting $R_{2h}u(t)$ from (\ref{R2hu}) into \added{the derivative on} the right-hand
side \changed{and} using the ODE system (\ref{ODEsystem2h}) to replace
$U'_{2h}(t)$, \added{we obtain
\begin{multline*}
\alpha_{2h}(t)=
\frac{2^q}{2^q-1}\Big(({U}_h^c)'(t)-F_{2h}(t,R_{2h}u(t))\Big)\\
+\frac1{2^q-1}\Big(F_{2h}(t,R_{2h}u(t))-F_{2h}(t,U_{2h}(t))\Big)+\gO(h^{q+1}).
\end{multline*}}%
\added{As \eqref{R2hu} and \eqref{UhcU2h} imply that
\[R_{2h}u(t)=U_h^c(t)-\frac1{2^q}\eta_{2h}(t)+\gO(h^{q+1})=U_{2h}(t)-\eta_{2h}(t)+\gO(h^{q+1})\]
}%
%Manipulating the expressions with
we get %
\begin{multline}
\alpha_{2h}(t) = \frac{2^q}{2^q-1}\Big((U_h^c)'(t)
-F_{2h}\Big(t,U_h^c(t)-\frac1{2^q}\eta_{2h}(t)+\gO(h^{q+1})\Big)\Big)\\
+\frac1{2^q-1}\Big(F_{2h}\Big(t,U_{2h}(t)-\eta_{2h}(t)+\gO(h^{q+1})\Big)
-F_{2h}(t,U_{2h}(t))\Big)+\gO(h^{q+1}). %
\end{multline}
Taylor expansions yield %
\bq \label{alphaexact}%
\alpha_{2h}(t) =
\frac{2^q}{2^q-1}\Big((U_h^c)'(t)-F_{2h}(t,U_h^c(t))\Big)+\gO(h^{q+1})\,. %
\eq %
Analogously to (\ref{globTimeError}), we set
$e_h^c(t)\!=\!V_h^c(t)-U_h^c(t)$. Substituting $(U_h^c)'(t)$ by $R^h_{2h}F_h(t,U_h(t))$ and
using again Taylor expansion it follows that
\begin{multline}
\alpha_{2h}(t) =
\frac{2^q}{2^q-1}\Big(R^h_{2h}F_h(t,V_h(t))-F_{2h}(t,V_h^c(t))\Big)+\gO(h^{q+1})\\
-\frac{2^q}{2^q-1}\Big(R^h_{2h}\big( \partial_{U_h}F_h(t,V_h(t))\,e_h(t)\big)
-\partial_{U_h}F_{2h}
(t,V_h^c(t))e_h^c(t)\Big)+\gO(e_h(t)^2)\,.
\end{multline} %
%For $q\!=\!2$, this formula was also given in \cite{Be1988}, Section
%6.1., (6.11).
Assuming the term on the right-hand
side involving the global time error to be sufficiently small, we can use%
\bq \label{alpha2happrox}%
\tilde{\alpha}_{2h}(t) =
\frac{2^q}{2^q-1}\Big(R^h_{2h}F_h(t,V_h(t))-F_{2h}(t,V_h^c(t))\Big) %
\eq %
as approximation for the spatial truncation error on the coarse mesh.
To guarantee a suitable quality of the estimate (\ref{alpha2happrox}) we shall
first control the global time error \changed{with the aim} that afterwards
the overall error is dominated by the spatial truncation error (see
Section \ref{secControlRule}).

An approximation $\tilde{\alpha}_{h}(t)$ of the spatial truncation
error on the (original) fine mesh is obtained by interpolation
respecting the order of accuracy (see Section \ref{secIntFormulas}).
Thus, to approximate the global
spatial error $\eta_h(t)$ we consider instead of
(\ref{odeSpatialErrLin}) the step-size frozen version %
\bq \label{odeSpaceErrFrozen} %
\begin{array}{rll} %
\teta'_h(t)&\!\!=\!\!&\partial_{U_h}F_h(t_n,V_{h,n})\,\teta_h(t)-
\tilde{\alpha}_h(t),\quad
t\in(t_n,t_{n+1}],\;n=0,\dots,M\!-\!1,\\[0.2cm]
\teta_h(0)&\!\!=\!\!&0\,.
\end{array} %
\eq %

\begin{Rem} \label{remarkonbetteralpha} \rm
If an approximation $\te_h(t)$ of the global time error has already
been computed, we could make use of $U_h^c(t)\approx
V_h^c(t)\changed{-}\te_h^c(t)$ to obtain a better approximation of
$\alpha_{2h}(t)$ from (\ref{alphaexact}). However, we have found \changed{the following in our
experiments: Using the step size frozen equations (\ref{odeTimeErrFrozen})
and (\ref{odeSpaceErrFrozen}) together with
(\ref{alphaexact}) to approximate the global time and
spatial error did not yield a significantly better
approximation, not even in the case when the global time error was not
small.} Since in practice the use of formula
(\ref{alphaexact}) requires additional function evaluations,
equation (\ref{alpha2happrox}) appears to be more efficient.
\hfill$\Diamond$ %
\end{Rem}

\begin{Rem} \label{remarkonboundary} \rm
We note that special care has to be taken in the handling the
spatial truncation error at the boundary when derivative boundary
conditions are present. This requests interpolation adopted to the
correct order of accuracy, see \cite{Be1988}.
\hfill$\Diamond$ %
\end{Rem}

\section{The example discretization formulas \label{secIntFormulas}}
In order to keep the illustration as simple as possible we restrict
ourselves to one space dimension. For the spatial discretization of
(\ref{PDEsystem}) we use standard second-order finite differences.
Hence we have $q\!=\!2$. The discrete $L_2$-norm on a non-uniform
mesh %
\bq\label{spatialMesh}
x_0<x_1<\ldots <x_N<x_{N+1}\,,\quad h_i=x_i-x_{i-1}\,,\quad i=1,\ldots,N+1\,,
\eq %
for a vector $y=(y_1,\ldots,y_N)^T\in\R^N$ is defined through %
\bq\label{discreteL2Norm}
\|y\|^2 = \sum_{i=1}^{N}\frac{h_i+h_{i+1}}{2}\,y_i^2\,.
\eq
Here, the components $y_0$ and $y_{N+1}$ which are given by the boundary
values are not considered.

\added{\it Adaptive time integration.} The example time integration formulas are taken from
\cite{LaVe2007}. For the sake of completeness we shall give a short
summary of the implementation used. To generate the time grid
(\ref{timegrid}) we use as an example integrator the
3rd-order, A-stable Runge-Kutta-Rosenbrock scheme ROS3P, see
\cite{LangBOOK,LaVe2001} for more details. The property of tolerance
proportionality \cite{ShampineBOOK} is
asymptotically ensured through working for the local residual with %
\bq \label{localTimeErrEst} %
Est = \frac23\,(I_h-\gamma\tau_nA_{h,n})^{-1}r_h(t_{n+1/2})\,,
\qquad A_{h,n}=\partial_{U_h}F_h(t_n,V_{h,n})\,, %
\eq %
where $\gamma$ is the stability coefficient of ROS3P. The common
filter $(I_h-\gamma\tau_nA_{h,n})$ serves to damp spurious stiff
components which would otherwise be amplified through the
$F_h$-evaluations within $r_h(t_{n+1/2})$.

Let $D_n\!=\!\|Est\|$ and $Tol_n\!=\!Tol_A+Tol_R\|V_{h,n}\|$ with
$Tol_A$ and $Tol_R$ given local tolerances. If $D_n>Tol_n$ the step
is rejected and
redone. Otherwise the step is accepted and we advance in time. In
both cases, \added{$r\tau_n$, where $r = (Tol_n/D_n)^{1/3}$, is in leading order equal to the step size which would have led to fulfill the local tolerance condition exactly, and which we therefore want to use in the next step. To be precautious, we multiply $r\tau_n$ with a safety factor of $0.9$. Further, to avoid too rapid step size changes, the step size is in each step only allowed to increase by maximally 50 \% and to decrease by maximally 1/3, leading overall to
the new step size being} determined by
\bq \label{tau-new}%
\tau_{new} = \mbox{min}\big(1.5,\mbox{max}(2/3, 0.9\, r)\big)\, %
\tau_n\,,\qquad r = (Tol_n/D_n)^{1/3}\,. %
\eq %
After each step size change we adjust $\tau_{new}$ to $\tau_{n+1} =
(T-t_n)/\lfloor(1+(T-t_n)/\tau_{new})\rfloor$ so as \changed{to avoid an
unnecessarily small final time step to reach the end point $T$.} The
initial step size $\tau_0$ is prescribed and is adjusted similarly.
\added{This heuristics works quite well in practice.}

The linear error transport equations (\ref{odeTimeErrFrozen}) and
(\ref{odeSpaceErrFrozen}) are simultaneously solved by means of the
implicit midpoint rule, which gives approximations $\te_{h,n}$ and
$\teta_{h,n}$ to the global time and spatial error at time
$t\!=\!t_n$. We use the implementations %
\bq \label{imr-time-error} %
\begin{array}{rll} %
(I_h - \frac12\tau_n A_{h,n})\, \delta e_{n+1} &=& 2\tilde{e}_{h,n} +
\frac23\tau_n r(t_{n+1/2})\,,%
\\ \rule{0cm}{4mm}
\tilde{e}_{h,n+1} &=& \delta e_{n+1} - \tilde{e}_{h,n}\,,
\end{array}
\eq %
and %
\bq \label{imr-space-error} %
\begin{array}{rll} %
(I_h - \frac12\tau_n A_{h,n})\, \delta\eta_{n+1} &=& 2\tilde{\eta}_{h,n} -
\tau_n \tilde{\alpha}_{h}(t_{n+1/2})\,,%
\\ \rule{0cm}{4mm}
\tilde{\eta}_{h,n+1} &=& \delta\eta_{n+1} - \tilde{\eta}_{h,n}\,.
\end{array}
\eq %
Clearly, the matrices $A_{h,n}$ already
computed within ROS3P can be reused. The spatial truncation error
$\tilde{\alpha}_{2h}(t)$ at $t\!=\!t_{n+1/2}$ is given by %
\bq %
\tilde{\alpha}_{2h}(t_{n+1/2})= %
\frac43\,\left(R_{2h}^hF_h\left(t_{n+1/2},V_h(t_{n+1/2})\right)- %
F_{2h}\left(t_{n+1/2},R_{2h}^hV_h(t_{n+1/2})\right)\right)\,.%
\eq %
Since $V_h(t_{n+1/2})$ and $F_h(t_{n+1/2},V_h(t_{n+1/2}))$ are
available from the computation of $r_h(t_{n+1/2})$ in
(\ref{reshalfwaydetailed}), this requires only one function evaluation
on the coarse grid. The vector
$\tilde{\alpha}_{2h}(t_{n+1/2})$ on the coarse mesh is prolongated
to the fine mesh and is then divided by $2^q\!=\!4$ if the neighbouring
fine grid points are equidistant, otherwise it is divided by $2^{q-1}\!=\!2$.
The remaining $\tilde{\alpha}_{h}(t_{n+1/2})$ on the fine mesh are
computed by interpolation respecting the order of the neighbouring
spatial truncation errors.

Due to freezing the coefficients in each time step, the second-order
midpoint rule is a first-order method when interpreted for solving
the linearised equations (\ref{odeTimeErrLin}) and (\ref{odeSpatialErrLin}). Thus if
all is going well, we asymptotically have
$\te_{h,n}\!=\!e_h(t_n)+\gO(\tau_{max}^4)$ and
$\teta_{h,n}\!=\!\eta_h(t_n)+\gO(\tau_{max}h_{max}^q)+\gO(h_{max}^{q+1})$.

After computing the spatial truncation errors we can solve the discretized
error transport equations (\ref{imr-space-error}) for all $\teta_{h,n}$. We
shall distinguish between two different mesh adaptation approaches: (i)
globally uniform and (ii) locally adaptive refinement. Although the
uniform strategy may be less efficient, it is very easy to implement and
therefore of special practical interest if software packages which do no allow
dynamic adaptive mesh refinement are used.

{\it Uniform spatial refinement.} Let $Tol$ be a given tolerance. Then our aim is to
guarantee $\|\eta_h(T)\|\le Tol$. From (\ref{imr-space-error}), we get an approximate
value $\teta_{h,M}$ for the spatial discretization error at $T$. If the
desired accuracy is still not satisfied, i.e., $\|\teta_{h,M}\|>Tol$, we choose a new
(uniform) spatial resolution
\bq \label{controlNewH}%
h_{new}=\sqrt[q]{\frac{Tol}{\|\teta_{h,M}\|}}\,h %
\eq %
to account for achieving $\|\eta_{h_{new}}(T)\|\approx Tol$. From $h_{new}$ we determine
a new number of mesh points. The whole computation is redone with the new
spatial mesh.

{\it Adaptive spatial refinement.} The main idea of our local spatial mesh control
is based on the observation that the principle of tolerance proportionality can
also be applied to the spatial discretization error. Multiplying all
$\tilde{\alpha}_{h}(t_{n+1/2})$ in (\ref{imr-space-error}) by a certain constant
multiplies all $\teta_{h,n+1}$ by the same constant since $\teta_{h,0}\!=\!0$.
Set $Tol_n^{\alpha}\!=\!Tol_A^{\alpha}+Tol_R^{\alpha}\|V_{h,n}\|$ where $Tol_A^{\alpha}$
and $Tol_R^{\alpha}$ are given local tolerances and define a local estimator $A_n$ through
\bq \label{spatialErrEst}%
A_n^2 = \sum_{i:\,x_i\in F_h} 2h_i|\tilde{\alpha}_i(t_{n+1/2})|^2,
\eq %
where $F_h$ denotes the set of all (fine) mesh points that do not belong to the coarse
mesh. Remember we have second order of the spatial truncation error in these points.
If $A_n\le Tol_n^{\alpha}$ the mesh is only coarsened. Otherwise, if $A_n>Tol_n^{\alpha}$
the mesh is improved by refinement and coarsening as well. We set
$\alpha_{tol}\!=\!0.9\,Tol_n^{\alpha}/\sqrt{N}$ and mark all $x_i\in F_h$
\bq \label{ref-strategy} %
\begin{array}{rll} %
\mbox{for refinement if } && \sqrt{h_i}\,\tilde{\alpha}_i(t_{n+1/2})>\alpha_{tol}%
\\ \rule{0cm}{4mm}
\mbox{and for coarsening if } && \sqrt{h_i}\,\tilde{\alpha}_i(t_{n+1/2})<0.1\,\alpha_{tol}\,.
\end{array}
\eq %
Grid adaptation is first performed for the coarse mesh and afterwards the fine mesh is
constructed by halving each interval. If $x_i$ is marked for refinement the corresponding
coarse grid interval is halved. Grid points are only removed if there are two equidistant
neighbouring intervals the midpoints of which are marked for coarsening. Finally, the grid
is smoothed such that $0.5\!\le\!h_i/h_{i-1}\!\le\!2$ everywhere. Data transfer from old to new
meshes is done by cubic Hermite interpolation where the necessary first derivatives are
determined from fourth order finite differences.

After mesh adaptation the local time step
is redone with the new mesh. The procedure is continued until first $D_n\le Tol_n$ and
second $A_n\le Tol_n^\alpha$ hold.
The whole strategy aims at equidistributing the local values
$\sqrt{h_i}\,\tilde{\alpha}_i(t_{n+1/2})$. Asymptotically we get
\bq \label{ref-strategy-asym} %
A_n \approx \left( 2 \sum_{i:\,x_i\in F_h} \alpha_{tol}^2 \right) ^{1/2} =
\left( 2 \sum_{i:\,x_i\in F_h} \frac{0.81\,(Tol_n^\alpha)^2}{N} \right) ^{1/2}
\approx 0.9\,Tol_n^{\alpha}\,,
\eq %
where the factor $0.9$ improves the robustness of the equidistribution principle.

\section{The control rules \label{secControlRule}}
Like for the ODE case studied in \cite{LaVe2007} our aim is to
provide global error estimates and to control the accuracy of the
numerically computed solution to the imposed tolerance level. Let
$GTol_A$ and $GTol_R$ be the global tolerances. Then we start with
the local tolerances $Tol_A=GTol_A$, $Tol_R=GTol_R$, and in the
spatially adaptive case also with
$Tol_A^\alpha=C_\alpha\,GTol_A$, and $Tol_R^\alpha=C_\alpha\,GTol_R$,
where the factor $C_\alpha>1$ ensures that the residual time error is
small with respect to the spatial truncation error and therefore the use of
(\ref{alpha2happrox}) is justified.

Suppose the numerical schemes have delivered an approximate solution
$V_{h,M}$ and global estimates $\te_{h,M}$ and $\teta_{h,M}$ for the
time and spatial error at time
$t_M\!=\!T$. We then verify whether %
\bq \label{controlRuleTime} %
\|\te_{h,M}\|\leq C_TC_{control}Tol_M,\quad
Tol_M=GTol_A+GTol_R\|V_{h,M}\|, %
\eq %
where $C_{control}\approx1$, typically $>1$, and $C_T\in(0,1)$
denotes the fraction desired for the global time error with respect
to the tolerance $Tol_M$. If (\ref{controlRuleTime}) does not hold,
the whole computation is redone over $[0,T]$ with the same initial
step $\tau_0$ and the adjusted local tolerances %
\bq \label{controlAdjTol} %
Tol_A=Tol_A\cdot fac,\quad Tol_R=Tol_R\cdot fac,\quad
fac=C_TTol_M/\|\te_{h,M}\|. %
\eq %
Based on tolerance proportionality, reducing the local error
estimates with the factor {\it fac} will reduce $e_h(T)$ by {\it
fac} \cite{ShampineBOOK}.

\begin{table}[ht]
\begin{center}
\setlength{\tabcolsep}{2mm}
\begin{tabular}{|l|l|}\hline
Step   & Control Algorithm with Uniform Refinement in Space
\\ \hline \hline%
Step 0 & Choose global tolerances $GTol_A$ and $GTol_R$.\\[0.5mm]
       & Choose $C_T$, $C_{control}$, $h_0$, $q$, and $\tau_0$.\\[0.5mm]
       & Set local tolerances $Tol_A=GTol_A$ and $Tol_R=GTol_R$.\\[0.5mm]
       & Set $h=h_0$.
\\ \hline \hline%
Step 1 & Run numerical schemes to compute $V_{h,M}, \te_{h,M}, \teta_{h,M}$.\\[0.5mm]
       & Compute $Tol_M=GTol_A+GTol_R\|V_{h,M}\|$.
\\ \hline \hline%
Step 2 & IF $\|\te_{h,M}\|\le C_TC_{control}Tol_M$ GOTO Step 3.\\[0.25mm]
       & ELSE set \\
       &\hspace{0.1cm} $fac=C_TTol_M/\|\te_{h,M}\|$, $Tol_A=Tol_A\cdot fac$, $Tol_R=Tol_R\cdot fac$\\[0.5mm]
       & and GOTO Step 1.
\\ \hline \hline%
Step 3 & IF $\|\te_{h,M}+\teta_{h,M}\|\le C_{control}Tol_M$ GOTO Step 4.\\[0.75mm]
       & ELSE set $h=\sqrt[q]{(1-C_T)Tol_M/\|\teta_{h,M}\|}\,h$ and GOTO Step 1.
\\ \hline \hline%
Step 4 & IF $h\ne h_0$ compute $q_{num}$.\\[0.5mm]
       & ELSE set $h=2h$, run numerical schemes again and compute then\\
       &\hspace{0.1cm} $q_{num}$.\\[0.5mm]
       & IF $q_{num}\approx q$ accept fine grid solution and STOP.\\[0.5mm]
       & ELSE set $h_0=2h_0, h=h_0$ and GOTO Step 1.
\\ \hline
\end{tabular}
\end{center}
\vspace{0.2cm}
\caption{Algorithmic structure of the overall control strategy when uniform refinement in
space is used.%
\label{table-control-algorithm-uniform}}
\end{table}

\begin{table}[ht]
\begin{center}
\setlength{\tabcolsep}{2mm}
\begin{tabular}{|l|l|}\hline
Step   & Control Algorithm with Adaptive Refinement in Space
\\ \hline \hline%
Step 0 & Choose global tolerances $GTol_A$ and $GTol_R$.\\[0.5mm]
       & Choose $C_T$, $C_{control}$, $C_\alpha$, $q$, and $\tau_0$.\\[0.5mm]
       & Set local tolerances\\
       & \hspace{0.1cm} $Tol_A\!=\!GTol_A$, $Tol_R\!=\!GTol_R$,
         $Tol_A^\alpha\!=\!C_\alpha\,GTol_A$, and\\
       & \hspace{0.1cm} $Tol_R^\alpha\!=\!C_\alpha\,GTol_R$.\\[0.5mm]
       & Choose initial spatial mesh.
\\ \hline \hline%
Step 1 & Run numerical schemes to compute $V_{h,M}, \te_{h,M}, \teta_{h,M}$.\\[0.5mm]
       & Compute $Tol_M=GTol_A+GTol_R\|V_{h,M}\|$.
\\ \hline \hline%
Step 2 & IF $\|\te_{h,M}\|\le C_TC_{control}Tol_M$ GOTO Step 3.\\[0.25mm]
       & ELSE set \\
       &\hspace{0.1cm} $fac=C_TTol_M/\|\te_{h,M}\|$, $Tol_A=Tol_A\cdot fac$, $Tol_R=Tol_R\cdot fac$\\[0.5mm]
       & and GOTO Step 1.
\\ \hline \hline%
Step 3 & IF $\|\te_{h,M}+\teta_{h,M}\|\le C_{control}Tol_M$ accept solution and STOP.\\[0.75mm]
       & ELSE set \\
       & \hspace{0.1cm} $fac\!=\!(1-C_T)Tol_M/\|\teta_{h,M}\|$, $Tol_A^\alpha\!=\!Tol_A^\alpha\cdot fac$,
         $Tol_R^\alpha\!=\!Tol_R^\alpha\cdot fac$\\[0.5mm]
       & and GOTO Step 1.
\\ \hline%
\end{tabular}
\end{center}
\vspace{0.2cm}
\caption{Algorithmic structure of the overall control strategy when adaptive refinement in
space is used.%
\label{table-control-algorithm-adaptive}}
\end{table}

If (\ref{controlRuleTime}) holds, we check whether %
\bq \label{controlRuleErr} %
\|\te_{h,M}+\teta_{h,M}\|\leq C_{control}Tol_M. %
\eq %
If it is true, the overall error
$E_h(T)\!=\!V_h(T)\!-\!(R_hu)(t)\!=\!e_h(T)\!+\!\eta_h(T)$ is
considered small enough relative to the chosen tolerance and
$V_{h,M}$ is accepted. Otherwise, the whole computation is redone
with the (already) adjusted tolerances (\ref{controlAdjTol}) and
an improved spatial resolution.

In the {\it uniform} case, we use the new mesh size computed from
(\ref{controlNewH}) with $Tol=(1-C_T)Tol_M$. To
check the convergence behaviour in space and therefore also the
quality of the approximation of the spatial truncation error, we
additionally compute the numerically observed order %
\bq \label{numObservedOrder} %
q_{num}=\log\left(\frac{\|\teta_{h,M}\|}{\|\teta_{h_{new},M}\|}\right)
\big/ \log\left(\frac{h}{h_{new}}\right).
\eq %
If $q_{num}$ computed for the final run is not close to the expected
value $q$ used for our Richardson extrapolation, we reason that the
approximation of the spatial truncation errors has failed due to a
dominating global time error, which happens, e.g., if the initial
spatial mesh is already too fine. Consequently, we coarsen the
initial mesh by a factor two and start again. If the control
approach stops without a mesh refinement, we perform an additional control
run on the coarse mesh and compute $q_{num}$ from
(\ref{numObservedOrder}) with $h_{new}\!=\!2h$. It turns out that
this simple strategy works quite robustly, provided that the meshes
used are able to resolve the basic behaviour of the solution.
The algorithmic structure of our control strategy with uniform refinement
in space is given in Table \ref{table-control-algorithm-uniform}.

In the {\it adaptive} case, we choose new local tolerances
\bq \label{controlAdjTolSpace} %
Tol_A^\alpha=Tol_A^\alpha\cdot fac,\quad Tol_R^\alpha=Tol_R^\alpha\cdot fac,\quad
fac=(1-C_T)Tol_M/\|\teta_{h,M}\|\,, %
\eq %
and the whole computation is redone over the interval $[0,T]$. Based on tolerance
proportionality, reducing the local truncation error with the factor {\it fac} will
reduce $\eta_h(T)$ by {\it fac}. In Table
\ref{table-control-algorithm-adaptive},
the algorithmic structure of our control strategy with adaptive refinement
in space is displayed. Note that now the index $h$ refers to a sequence of
spatial meshes adapted at each time point $t_n$.

Summarizing, the first check (\ref{controlRuleTime}) and the
possible second control computation serve to significantly reduce
the global time error. This enables us to make use of the
approximation (\ref{alpha2happrox}) for the spatial truncation
error, which otherwise could not be trusted. The second step based on
suitable spatial mesh improvement attempts to bring the overall
error down to the imposed tolerance. Using the sum of the approximate global
time and spatial error inside the norm in (\ref{controlRuleErr}), we
take advantage of favourable effects of error cancellation. These two
steps are successively repeated until the second check is
successful. Additionally, if uniform mesh refinement is used we take
into account the numerically observed order in space to assess the
approximation of the spatial truncation error.

\section{Numerical illustrations \label{secNumExamples}}
To illustrate the performance of the global error estimators and the
control strategy, we consider three test problems:
(i) the highly stable heat equation with
nonhomogeneous Neumann boundary conditions~\cite{Be1988}, (ii)
the nonlinear convection-dominated Burgers'
equation~\cite{Be1988,LaBeDe1991}, and (iii) the
Allen-Cahn equation modelling a diffusion-reaction
problem~\cite{LaVe2007}. Analytic solutions are known for
all three problems. Uniform spatial refinement is studied for
all three test cases. For \added{the Burgers' and} Allen-Cahn problem, these
results are compared to those obtained with adaptive refinement.
\added{We omit the corresponding results for the heat equation, since
the solution is very smooth in space and hence adaptive refinement
is not necessary. The challenge here is to control
the fast decay in time.}

We set $GTol_A=GTol_R=GTol$ for $GTol=10^{-l},l=2,\ldots,7$ and
start with one and the same initial step size $\tau_0=10^{-5}$.
Equally spaced meshes of $25$, $51$, $103$, $207$, $415$, $831$, and
$1663$ points are used as initial mesh. The control parameters
introduced above for the control rules are $C_T=1/3$,
$C_{control}=1.2$, and $C_\alpha=10$. All runs were performed, but
for convenience we only select a representative set of them for our
presentation.

We define the estimated global error
$\tilde{E}_{h,M}=\te_{h,M}+\teta_{h,M}$ at time $t=T$ and set
indicators $\Theta_{est}=\|\tilde{E}_{h,M}\|/\|E_h(T)\|$ for the
ratio of the estimated global error and the true global error, and
$\Theta_{ctr}=Tol_M/\|E_h(T)\|$ for the ratio of the desired
tolerance and the true global error. Thus, $\Theta_{ctr}\ge
1/C_{control}=5/6$ indicates control of the true global error.

The tables of results contain the following quantities,
$Tol=Tol_A=Tol_R$ from (\ref{controlAdjTol}),
$Tol^\alpha=Tol_A^\alpha=Tol_R^\alpha$ from (\ref{controlAdjTolSpace}),
$Tol_M = GTol\,(1 + \|V_{h,M}\|)$ from (\ref{controlRuleTime}), the estimated
global error $\tilde{E}_{h,M}$, the estimated time error
$\te_{h,M}$, and the estimated spatial truncation error
$\teta_{h,M}$. Note that we always start with $Tol=GTol$ in the first run.
The ratios $\Theta_{est}$ and $\Theta_{ctr}$ serve to
illustrate the quality of the global error estimation and the
control. If uniform refinement in space is applied, the numerically observed
order $q_{num}$ for the spatial error is given. It will be clear from the
tables of results whether a tolerance-adapted run to control the
global time error, a spatial mesh adaptation step or an additional
control run on a coarser grid was necessary. Especially, the latter
is marked by a dashed line.

\begin{table}[ht]
\begin{center}
\setlength{\tabcolsep}{1.99mm}
\begin{tabular}{|r|r||c|c|c|c||r|r|r|}\hline
$Tol$  & \rule{0mm}{0.4cm}$N$ & $Tol_M$ & $\|\tilde{E}_{h,M}\|$ &
$\|\te_{h,M}\|$ & $\|\teta_{h,M}\|$ & $\Theta_{est}$ &
$\Theta_{ctr}$ & $q_{num}$
\\ \hline \hline%
1.00e-2&  25&1.10e-2&7.14e-4&1.16e-4&8.20e-4&  0.99& 15.27&
\\\hdashline1.00e-2&  13&1.10e-2&3.27e-3&1.24e-4&3.38e-3& 0.99&  3.32&  2.04
\\ \hline\hline%
1.00e-3&  51&1.10e-3&1.68e-4&1.97e-5&1.86e-4&  1.00&  6.51&
\\\hdashline1.00e-3&  25&1.10e-3&8.04e-4&2.03e-5&8.22e-4& 1.00&  1.36&  2.02\\ \hline\hline%
1.00e-4& 103&1.10e-4&4.27e-5&2.01e-6&4.44e-5& 1.00&  2.57&
\\\hdashline1.00e-4&  51&1.10e-4&1.85e-4&1.96e-6&1.86e-4&1.00&  0.59&  2.01
\\ \hline\hline%
1.00e-5& 207&1.10e-5&1.07e-5&1.89e-7&1.08e-5&1.00&  1.03&
\\\hdashline1.00e-5& 103&1.10e-5&4.43e-5&1.83e-7&4.44e-5&1.00&  0.25&  2.01
\\ \hline\hline%
1.00e-6& 415&1.10e-6&2.67e-6&1.81e-8&2.68e-6&  1.00&  0.41&
\\1.00e-6& 795&1.10e-6&7.14e-7&1.83e-8&7.28e-7&1.00&  1.54&  2.00
\\ \hline\hline%
1.00e-7&  25&1.10e-7&8.24e-4&1.24e-9&8.24e-4& 1.00&  0.00&
\\1.00e-7&2759&1.10e-7&5.91e-8&1.60e-9&6.03e-8&1.00&  1.86&  2.01
\\ \hline\hline%
1.00e-7&1663&1.10e-7&1.65e-7&1.57e-9&1.66e-7& 1.00&  0.67&
\\1.00e-7&2505&1.10e-7&7.20e-8&1.57e-9&7.31e-8&1.00&  1.53&  2.00\\ \hline
\end{tabular}
\end{center}
\vspace{0.2cm}
\caption{Selected data for the heat equation with Neumann boundary conditions.
Uniform refinement in space is used.%
\label{table-test-problem2}}
\end{table}

\subsection{Heat equation with Neumann boundary conditions}
This heat equation provides an example with inhomogeneous Neumann
boundary conditions: %
\bq \label{test-problem2} %
\partial_tu = \partial_{xx}u\,, \quad 0 < x < 1.0\,,
\qquad 0 < t \leq T = 0.2\,,
\eq %
and boundary conditions $\partial_xu=\pi\,e^{-\pi^2t}\cos(\pi x)$ at
$x\!=\!0$ and $x\!=\!1$. The initial condition is consistent with
the analytic solution $u(x,t)\!=\!e^{-\pi^2t}\sin(\pi x)$. Although
the solution is very stable, it is not easy to provide good error
estimates as stated in \cite{Be1988,LaBeDe1991}.

To approximate the inhomogeneous Neumann boundary conditions we introduce
artificial mesh points $x_{-1}=-h$ and $x_{N+2}=1+h$, discretize
$\partial_xu(0)$ and $\partial_xu(1)$ by second
order central differences, and use the approximate differential equation at
the boundary to eliminate the artificial solution values.
In consequence, we have global spatial order $q=2$ in all mesh points, but when interpolating the estimated spatial truncation error we have to respect that it is of first order at the boundary (see also Remark \ref{remarkonboundary}).

Due to the high stability of the problem the global time errors are
much smaller than imposed local tolerances. So, control of the
global time error is redundant here and control runs were only
carried out in case of insufficient spatial resolutions.
Table~\ref{table-test-problem2} shows results for various tolerances
and initial meshes. \added{We select two runs to explain the control strategy.
For the third simulation, we take $GTol=10^{-4}$ and start with the local
tolerance $Tol=10^{-4}$. Using $103$ mesh points in space, we run the
computation and get the following approximations of the time and spatial
errors: $\|\te_{h,M}\|=2.01\times 10^{-6}$ and
$\|\teta_{h,M}\|=4.44\times 10^{-5}$. The control checks for the time error
estimate, $\|\te_{h,M}\|\le C_TC_{control}Tol_M=4.4\times 10^{-5}$,
and for the global error,
$\|\tilde{E}_{h,M}\|=4.27\times 10^{-5}\le 1.32\times 10^{-4}=C_{control}Tol_M$,
are positive, so that we already can stop after the first run. In accordance
with our safety strategy, we additionally perform one run on a coarser
mesh with half of the grid points, i.e., $N=51$. The numerically observed
order computed from (\ref{numObservedOrder}) is $q_{num}=2.01$. We reason that
our assumption for a successful Richardson extrapolation to estimate
the spatial truncation error is fulfilled and accept the numerical solution.
Choosing $GTol=10^{-7}$ and $N=25$, the approximate time error is still
very small, but the check for the global error,
$8.24\times 10^{-4}\le 1.32\times 10^{-7}$, obviously fails. From (\ref{controlNewH}),
we compute a new number of spatial mesh points, $N=2759$. Finally, the second
run is successful and with the numerically observed spatial order $q_{num}=2.01$ the
numerical solution is accepted.}

\begin{figure}[ht]
\centering
\includegraphics[width=0.40\textwidth]{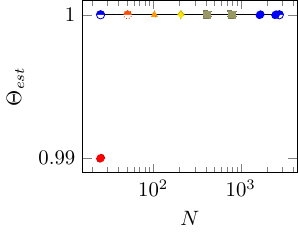}
\hspace{0.35cm}
\includegraphics[width=0.52\textwidth]{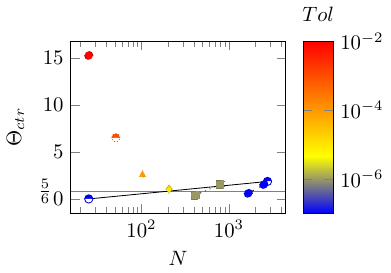}
\caption{Evolution of the efficiency indicators $\Theta_{est}$ (left) and $\Theta_{ctr}$
(right) for the heat equation with Neumann boundary conditions and
global tolerances $GTol=10^{-i},\, i=2,\ldots,7$. Different icons represent
different $GTols$. The progress in the local time tolerance $Tol$ is described
by diverse colouring. Control of the true global error,
i.e.\ $\Theta_{ctr}\ge 5/6$, is achieved in all cases.
Only for higher tolerances $GTol=10^{-6},\, 10^{-7}$, a second run is necessary,
indicated by connected icons. The quality of the estimates is very high.}
\label{fig-effHeat}
\end{figure}

The global error estimation and control appear
to work very well for this problem, where the influence of the
initial mesh points is less strong. This holds also for other
combinations of tolerances and initial meshes. \added{The results
are visualized in Fig.~\ref{fig-effHeat}}. Note the high quality
of the estimator $\tilde{E}_{h,M}$ (and therefore also of $\teta_{h,M}$),
showing that the derivative boundary
condition is well resolved within the Richardson extrapolation. For
the runs with tolerances $GTol\!=\!10^{-2}, 10^{-3},1 0^{-4},
10^{-5}$, the order of the spatial convergence was successfully
checked with a second run on the coarse mesh, that is, we can trust
the first run.

\begin{table}[ht]
\begin{center}
\setlength{\tabcolsep}{2mm}
\begin{tabular}{|r|r||c|c|c|c||r|r|r|}\hline
$Tol$  & \rule{0mm}{0.4cm}$N$ & $Tol_M$ & $\|\tilde{E}_{h,M}\|$ &
$\|\te_{h,M}\|$ & $\|\teta_{h,M}\|$ & $\Theta_{est}$ &
$\Theta_{ctr}$ & $q_{num}$
\\ \hline \hline%
1.00e-2&  51&1.93e-2&4.30e-3&1.86e-3&2.86e-3&  1.08&  4.87&
\\\hdashline1.00e-2&  25&1.93e-2&1.29e-2&2.21e-3&1.14e-2&0.99&  1.48&  2.00\\ \hline\hline%
1.00e-3&  51&1.93e-3&2.83e-3&1.54e-4&2.74e-3& 0.99&  0.68&
\\1.00e-3&  75&1.93e-3&1.36e-3&1.48e-4&1.28e-3&1.00&  1.42&  2.00\\ \hline\hline%
1.00e-4&  51&1.93e-4&2.73e-3&1.09e-5&2.73e-3&  0.98&  0.07&
\\1.00e-4& 239&1.94e-4&1.32e-4&1.05e-5&1.27e-4&1.00&  1.46&  2.00\\ \hline\hline%
1.00e-5&  51&1.93e-5&2.73e-3&1.08e-6&2.73e-3&  0.98&  0.01&
\\1.00e-5& 757&1.94e-5&1.32e-5&1.02e-6&1.27e-5&1.00&  1.47&  2.00\\ \hline\hline%
1.00e-6&  51&1.93e-6&2.73e-3&1.08e-7&2.73e-3& 0.98&  0.00&
\\1.00e-6&2391&1.94e-6&1.32e-6&9.29e-8&1.28e-6&1.00&  1.47&  2.00\\ \hline\hline%
1.00e-7&  51&1.93e-7&2.73e-3&1.10e-8&2.73e-3&0.98&  0.00&
\\1.00e-7&7563&1.94e-7&1.31e-7&8.57e-9&1.28e-7&1.00&  1.47&  2.00
\\ \hline
\end{tabular}
\end{center}
\vspace{0.2cm}
\caption{Selected data for Burgers' equation with $51$ initial mesh points.
Uniform refinement in space is used.%
\label{table-test-problem3a}}
\end{table}

\begin{table}[ht]
\begin{center}
\setlength{\tabcolsep}{1.9mm}
\begin{tabular}{|r|r|r||c|c|c|c||r|r|}\hline
$Tol$  & $Tol^\alpha$  &\rule{0mm}{0.4cm}$N_M$ & $Tol_M$ & $\|\tilde{E}_{h,M}\|$ &
$\|\te_{h,M}\|$ & $\|\teta_{h,M}\|$ & $\Theta_{est}$ &
$\Theta_{ctr}$
\\ \hline \hline%
1.00e-2&1.00e-1&  15&1.92e-2&2.81e-2&3.15e-3&2.61e-2&  1.64&  1.12\\
1.00e-2&4.91e-2&  25&1.93e-2&1.57e-2&2.30e-3&1.46e-2&  1.17&  1.44\\
\hline\hline%
1.00e-3&1.00e-2&  45&1.92e-3&9.93e-4&1.15e-4&9.46e-4&  1.01&  1.95\\
\hline\hline%
1.00e-3&1.00e-1&  13&1.90e-3&1.49e-2&2.29e-4&1.48e-2&  1.04&  0.13\\
1.00e-3&8.53e-3&  49&1.92e-3&8.18e-4&1.14e-4&7.85e-4&  1.03&  2.41\\
\hline\hline%
1.00e-4&1.00e-2&  43&1.92e-4&7.95e-4&1.03e-5&7.93e-4&  1.02&  0.25\\
1.00e-4&1.61e-3&  89&1.92e-4&1.94e-4&1.09e-5&1.89e-4&  1.01&  1.00\\
\hline
\end{tabular}
\end{center}
\vspace{0.2cm}
\caption{Selected data for Burgers' equation. Adaptive refinement
in space is used.%
\label{table-test-problemBurgers-adaptive}}
\end{table}

\subsection{Burgers' equation}
The second problem is the nonlinear Burgers' equation %
\bq \label{test-problem3} %
\partial_tu = \varepsilon\,\partial_{xx}u - u\partial_xu\,,\quad 0 < x < 1.0\,,
\qquad 0 < t \leq T = 1.0\,, %
\eq %
where $\varepsilon=0.015$ is used in the experiments. Dirichlet
boundary conditions and initial conditions are consistent with the analytic
solution defined by %
\bq %
u(x,t) = \frac{r_1+5r_2+10r_3}{10(r_1+r_2+r_3)}\,,
\eq %
where $r_1(x)=e^{0.45x/\varepsilon}$,
$r_2(t,x)=e^{0.01(10+6t+25x)/\varepsilon}$, and
$r_3(t)=e^{0.025(6.5+9.9t)/\varepsilon}$.

We note that this equation does not formally fit into our setting of semilinear parabolic equations\added{, and e.g.\ the linearized error transport equations \eqref{odeSpatialErrLin} and \eqref{odeTimeErrLin} are no longer valid, as the $\gO$-terms would now be divided by the spatial discretization step size $h$}. However, it is indeed interesting to see how the proposed algorithm performs for this widely used benchmark problem.

In Table~\ref{table-test-problem3a} we present results \added{with uniform refinement
in space} for all tolerances used and the
$51$-point initial mesh. The use of a relatively coarse mesh at the beginning is the
natural choice in practice. No adaptation in time is necessary, which is mainly due to
the small first time step and the maximum factor $1.5$ which is allowed in (\ref{tau-new})
for a step size enlargement. For the tolerance $GTol\!=\!10^{-2}$, the numerical
solution is accepted since the corresponding control run on a coarser mesh shows $q_{num}\!\approx\!2$,
the expected value. Remarkably excellent estimates are obtained for higher tolerances.
Here, control is always achieved after one spatial mesh improvement.

\added{Let us have a closer look at the second run. We choose
$GTol=10^{-3}$ and start with a local tolerance $Tol=10^{-3}$ for the time integrator. The
inspection of the global time error, $\|\te_{h,M}\|=1.54\times 10^{-4}$, shows that the
control rule (\ref{controlRuleTime}) is fulfilled. So, an adaption of the local tolerance $Tol$ is
not necessary. However, the approximate global error, $\|\tilde{E}_{h,M}\|=2.83\times 10^{-3}$, is
still too large due to an unacceptable spatial error, $\|\teta_{h,M}\|=2.74\times 10^{-3}$.
We compute a new number of spatial points, $N=75$, from (\ref{controlNewH}) and perform a
second run which is now successful. With the numerically observed spatial order $q_{num}=2.00$ the
numerical solution is considered as accurate enough.}

\begin{figure}[ht]
\centering
\includegraphics[width=0.415\textwidth]{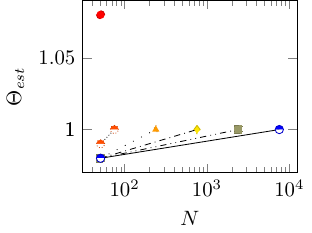}
\hspace{0.35cm}
\includegraphics[width=0.52\textwidth]{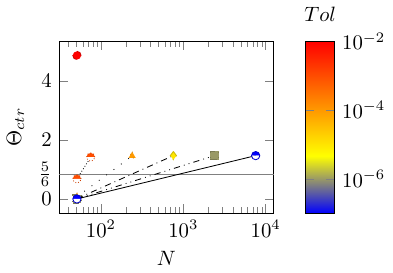}
\caption{Evolution of the efficiency indicators $\Theta_{est}$ (left) and $\Theta_{ctr}$
(right) for the Burgers problem, global tolerances $GTol=10^{-i},\, i=2,\ldots,7$,
and uniform refinement in space. Different icons represent
different $GTols$. The progress in the local time tolerance $Tol$ is described
by diverse colouring. Control of the true global error,
i.e.\ $\Theta_{ctr}\ge 5/6$, is achieved in all cases.
Except for $GTol=10^{-2}$, a second run is necessary for all global tolerances,
indicated by connected icons. The quality of the estimates is very high.}
\label{fig-effBU}
\end{figure}

 \added{The evolution of the indicators
$\Theta_{est}$ and $\Theta_{ctr}$ is visualized in Fig.~\ref{fig-effBU}.}

\begin{figure}[ht]
\centering
\includegraphics[width=0.39\textwidth]{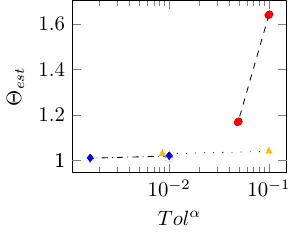}
\hspace{0.35cm}
\includegraphics[width=0.52\textwidth]{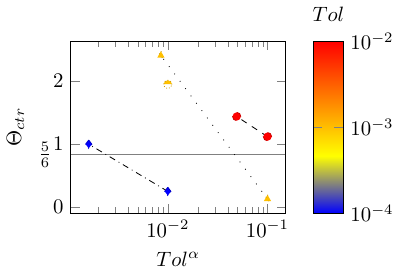}
\caption{Evolution of the efficiency indicators $\Theta_{est}$ (left) and $\Theta_{ctr}$
(right) for the Burgers problem, global tolerances $GTol=10^{-i},\, i=2,\ldots,4$,
and adaptive refinement in space. Here, $Tol^\alpha$ is the local spatial tolerance.
Different icons represent different $GTols$. The progress in the local time tolerance
$Tol$ is described by diverse colouring. Control of the true global error,
i.e.\ $\Theta_{ctr}\ge 5/6$, is achieved in all cases.
Except for $GTol=10^{-3}$ and $Tol^\alpha=10^{-2}$, a second run is necessary for all global tolerances,
indicated by connected icons. The quality of the estimates is very high.}
\label{fig-effBU-adapt}
\end{figure}

\added{The overall algorithm performs also well
when adaptive spatial refinement is used, as can be seen from Table~\ref{table-test-problemBurgers-adaptive}.
The quality of the estimation process is again very good, which leads to a significant reduction of
the number of mesh points compared with the uniform approach. We have used $C_\alpha=10$ in the first two
runs and $C_\alpha=100$ in the other ones
to set $Tol^\alpha=C_\alpha\,Tol$ at the beginning. The number of adaptive grid points at the final time $T$
is denoted by $N_M$. After adjusting the spatial meshes
until $A_n\le Tol_n^\alpha=Tol^\alpha(1+\|V_{h,n}\|)$ holds, no further runs with higher tolerances
in time are necessary. The evolution of the indicators
$\Theta_{est}$ and $\Theta_{ctr}$ is visualized in Fig.~\ref{fig-effBU-adapt}.}

\begin{figure}[ht]
\centering
\includegraphics[width=0.49\textwidth]{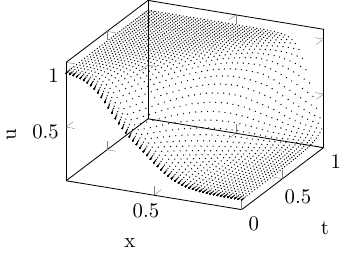}
\includegraphics[width=0.49\textwidth]{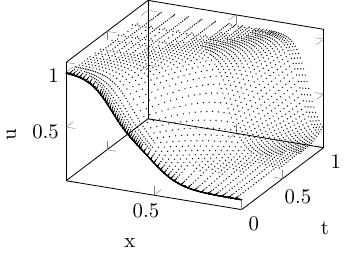}
\caption{Temporal evolution of the numerical solution for
the Burgers' problem with $Tol=10^{-3}$ and $51$ uniform grid
points (left) and adaptive spatial refinement with $45$
grid points at the final time (right).}
\label{fig-solmeshBU}
\end{figure}

\added{The numerical solutions obtained with $Tol=10^{-3}$ and $51$ uniform grid
points (left) and adaptive spatial refinement with $45$ grid points at the final time (right)
are plotted in Fig.~\ref{fig-solmeshBU}. With less grid points, the adaptive scheme reduces
the global error by nearly a factor $3$.}

\begin{table}[ht]
\begin{center}
\setlength{\tabcolsep}{2mm}
\begin{tabular}{|r|r||c|c|c|c||r|r|r|}\hline
$Tol$  & \rule{0mm}{0.4cm}$N$ & $Tol_M$ & $\|\tilde{E}_{h,M}\|$ &
$\|\te_{h,M}\|$ & $\|\teta_{h,M}\|$ & $\Theta_{est}$ &
$\Theta_{ctr}$ & $q_{num}$
\\ \hline \hline%
1.00e-2& 103&2.05e-2&1.84e-0&1.45e-1&1.98e-0&    9.89&  0.11&
\\4.69e-4& 103&2.05e-2&5.78e-1&1.26e-3&5.79e-1&  2.69&  0.10&
\\4.69e-4& 677&2.02e-2&6.04e-3&1.11e-3&7.15e-3&  1.19&  3.98&  2.34
\\ \hline\hline%
1.00e-2& 415&2.02e-2&7.69e-2&1.44e-1&6.73e-2&    3.05&  0.80&
\\4.66e-4& 415&2.02e-2&1.86e-2&1.11e-3&1.97e-2&  1.23&  1.34&
\\\hdashline4.66e-4& 207&2.03e-2&9.17e-2&1.15e-3&9.29e-2& 1.47&  0.32&  2.24
\\ \hline\hline%
1.00e-3& 207&2.03e-3&9.82e-2&2.97e-3&1.01e-1&  1.60&  0.03&
\\2.27e-4& 207&2.03e-3&8.80e-2&4.93e-4&8.85e-2&1.39&  0.03&
\\2.27e-4&1683&2.02e-3&6.14e-4&4.71e-4&1.09e-3&1.11&  3.67&  2.10
\\ \hline\hline%
1.00e-3& 831&2.02e-3&2.26e-3&2.87e-3&5.12e-3&  1.33&  1.19&
\\2.35e-4& 831&2.02e-3&4.01e-3&4.91e-4&4.50e-3& 1.12&  0.57&
\\2.35e-4&1521&2.02e-3&8.42e-4&4.90e-4&1.33e-3& 1.12&  2.68&  2.02
\\ \hline\hline%
1.00e-4&1663&2.02e-4&8.89e-4&1.86e-4&1.08e-3&  1.07&  0.24&
\\3.63e-5&1663&2.02e-4&9.88e-4&6.14e-5&1.05e-3&1.05&  0.21&
\\3.63e-5&4643&2.02e-4&7.30e-5&6.14e-5&1.34e-4&1.04&  2.89&  2.00
\\ \hline
\end{tabular}
\end{center}
\vspace{0.2cm}
\caption{Selected data for the Allen-Cahn problem. Uniform refinement
in space is used.%
\label{table-test-problem1}}
\end{table}

\begin{table}[ht]
\begin{center}
\setlength{\tabcolsep}{1.9mm}
\begin{tabular}{|r|r|r||c|c|c|c||r|r|}\hline
$Tol$  & $Tol^\alpha$  &\rule{0mm}{0.4cm}$N_M$ & $Tol_M$ & $\|\tilde{E}_{h,M}\|$ &
$\|\te_{h,M}\|$ & $\|\teta_{h,M}\|$ & $\Theta_{est}$ &
$\Theta_{ctr}$
\\ \hline \hline%
1.00e-2&1.00e-1&245&2.01e-2&1.05e-1&1.39e-1&3.42e-2&3.21&0.61\\
4.81e-4&1.00e-1&247&2.01e-2&8.54e-3&1.13e-3&9.67e-3&1.26&2.98\\
\hline\hline%
1.00e-3&1.00e-2&483&2.01e-3&1.26e-3&2.86e-3&1.59e-3&1.26&2.01\\
2.35e-4&1.00e-2&481&2.01e-3&9.72e-4&4.84e-4&1.46e-3&1.11&2.30\\
\hline\hline%
1.00e-4&1.00e-3&1839&2.01e-4&9.08e-5&1.85e-4&9.45e-5&1.63&3.62\\
3.62e-5&1.00e-3&1839&2.01e-4&5.49e-5&6.06e-5&1.16e-4&0.92&3.36\\
\hline\hline%
1.00e-4&1.00e-2& 481&2.01e-4&1.23e-3&1.85e-4&1.41e-3&1.08&0.18\\
3.62e-5&1.00e-2& 483&2.01e-4&1.32e-3&6.09e-5&1.38e-3&1.06&0.16\\
3.62e-5&9.68e-4&1839&2.01e-4&5.48e-5&6.07e-5&1.15e-4&0.92&3.36\\
\hline\hline%
1.00e-4&1.00e-1& 243&2.01e-4&8.66e-3&1.84e-4&8.84e-3&1.15&0.03\\
3.62e-5&1.00e-1& 243&2.01e-4&8.55e-3&6.08e-5&8.61e-3&1.12&0.03\\
3.62e-5&1.55e-3&1809&2.01e-4&5.68e-5&6.06e-5&1.17e-4&0.92&3.25\\
\hline%
\end{tabular}
\end{center}
\vspace{0.2cm}
\caption{Selected data for the Allen-Cahn problem. Adaptive refinement
in space is used.%
\label{table-test-problem1-adaptive}}
\end{table}

\subsection{The Allen-Cahn equation}
The third problem is the bi-stable Allen-Cahn equation which is
defined by
\bq \label{test-problem1} %
\partial_tu = 10^{-2}\,\partial_{xx}u + 100u \,(1-u^2)\,, \quad 0 < x <
2.5\,,\qquad 0 < t \leq T = 0.5\,,
\eq %
with the initial function and Dirichlet boundary values taken from
the exact wave front solution $u(x,t) = (1+e^{\lambda\,(x -
\alpha\,t)})^{-1}, \,\lambda = 50\,\sqrt{2}, \,\alpha =
1.5\,\sqrt{2}$. This problem was also used in \cite{DeLa2013,LaVe2007}.

First we apply uniform refinement in space.
Table \ref{table-test-problem1} reveals a high quality of the global
error estimation and also the control process works quite well. Let
us pick one exemplary run out to explain the overall control strategy
in more detail. Starting with $GTol=Tol=10^{-3}$ and $831$ mesh points,
which corresponds to the fourth simulation, the numerical scheme
delivers global error estimates $\|\te_{h,M}\|=2.87\times 10^{-3}$ and
$\|\teta_{h,M}\|=5.12\times 10^{-3}$ for the time and spatial error
of the approximate solution $V_{h,M}$ at the final time $t_M=T$. The
first check for the time error estimate
$\|\te_{h,M}\|\le C_TC_{control}Tol_M=8.08\times 10^{-4}$ fails
and we adjust the local tolerances by a factor
$fac=C_TTol_M/\|\te_{h,M}\|=2.35\times 10^{-1}$, which yields the new
$Tol=2.35\,10^{-4}$. The computation is then redone. Due to the tolerance
proportionality, in the second run the time error is significantly reduced
and the inequality $\|\te_{h,M}\|\le 8.08\times 10^{-4}$ is now
valid. We proceed with checking
$\|\tilde{E}_{h,M}\|\le C_{control}Tol_M=2.42\times 10^{-3}$, which is
still not true. From (\ref{controlNewH}), we compute a new number
of spatial mesh points $N=1521$. Finally, the third run is successful and
with the numerically observed spatial order $q_{num}=2.02$ the
numerical solution is accepted.

\begin{figure}[ht]
\centering
\includegraphics[width=0.39\textwidth]{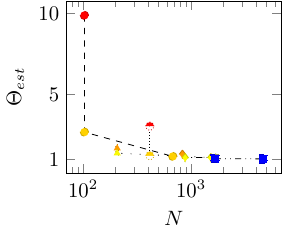}
\hspace{0.35cm}
\includegraphics[width=0.52\textwidth]{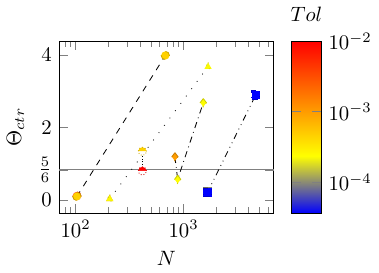}
\caption{Evolution of the efficiency indicators $\Theta_{est}$ (left) and $\Theta_{ctr}$
(right) for the Allen-Cahn problem, global tolerances $GTol=10^{-i},\, i=2,\ldots,4$,
and uniform refinement in space. Different icons represent
different $GTols$. The progress in the local time tolerance $Tol$ is described
by diverse colouring. Control of the true global error,
i.e.\ $\Theta_{ctr}\ge 5/6$, is achieved in all cases.
An improvement of the spatial meshes is necessary for all global tolerances,
indicated by connected icons. The quality of the estimates is very high after
the control runs.}
\label{fig-effAC}
\end{figure}

\begin{figure}[ht]
\centering
\includegraphics[width=0.39\textwidth]{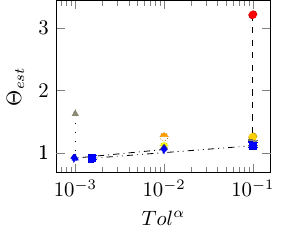}
\hspace{0.35cm}
\includegraphics[width=0.545\textwidth]{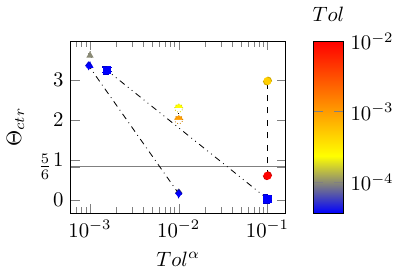}
\caption{Evolution of the efficiency indicators $\Theta_{est}$ (left) and $\Theta_{ctr}$
(right) for the Allen-Cahn problem, global tolerances $GTol=10^{-i},\, i=2,\ldots,4$,
and adaptive refinement in space. Here, $Tol^\alpha$ is the local spatial tolerance.
Different icons represent different $GTols$. The progress in the local time tolerance
$Tol$ is described by diverse colouring. Control of the true global error,
i.e.\ $\Theta_{ctr}\ge 5/6$, is achieved in all cases. The quality of the estimates is very high.}
\label{fig-effAC-adapt}
\end{figure}

The ratios for $\Theta_{est}=\|\tilde{E}_{h,M}\|/\|E_h(T)\|$
\added{in Table \ref{table-test-problem1}} lie between
1.04 and 1.23, after the control runs. Control of the global error,
that is $\|E_h(T)\|\leq C_{control}Tol_M$, is in general achieved
after two steps (one step to adjust the time grid and one step to
control the spatial discretization), whereas the efficiency index
$\Theta_{ctr}=Tol_M/\|E_h(T)\|$ is close to three. This results from
a systematic cancellation effect between the global time and spatial error,
which is not taken into account when computing $h_{new}$ from (\ref{controlNewH}).
\added{The evolution of the indicators
$\Theta_{est}$ and $\Theta_{ctr}$ is visualized in Fig.~\ref{fig-effAC}.}

\begin{figure}[ht]
\centering
\includegraphics[width=0.49\textwidth]{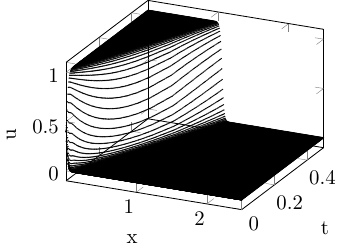}
\includegraphics[width=0.49\textwidth]{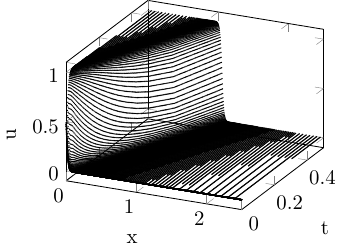}
\caption{Temporal evolution of the numerical solution for
the Allen-Cahn problem with $Tol=10^{-2}$ and $415$ uniform grid
points (left) and adaptive spatial refinement with $245$
grid points at the final time (right).}
\label{fig-solmeshAC}
\end{figure}

Next we consider locally adaptive spatial grid enhancement instead of globally
uniform adaptation. Within each time step the grid is adapted by refinement and
coarsening, based on an equidistribution principle, until
$A_n\le Tol_n^\alpha=Tol^\alpha(1+\|V_{h,n}\|)$ holds. This yields
a sequence of non-uniform meshes. Let $N_M$ denote the number of adaptive
grid points obtained at the final time $T$. The first three runs in
Table \ref{table-test-problem1-adaptive} correspond to our standard setting
$C_\alpha=10$, i.e.,
$Tol^\alpha=10\,Tol$. In this case, after adjusting the local tolerances
for the time integration no further run with higher tolerances in space
is necessary. To demonstrate the robustness of the algorithm, we select
two additional runs with $C_\alpha=10^l, l=2,3$, for $GTol=10^{-4}$. In both
cases, coarser meshes are used at the beginning and a second control run has to be
done to decrease the spatial discretization error. The resulting adaptive
spatial meshes are comparable. Control of the global error is always achieved.
The estimation process works again quite well. \added{The evolution of the indicators
$\Theta_{est}$ and $\Theta_{ctr}$ is visualized in Fig.~\ref{fig-effAC-adapt}.}
Compared to the uniform
case, significantly less spatial degrees of freedoms are needed to reach
the desired tolerances. The reduction rate varies between 40\% and 70\%.
\added{In Fig.~\ref{fig-solmeshAC} we have plotted the numerical solutions
obtained with $Tol=10^{-2}$ and $415$ uniform grid
points (left) and adaptive spatial refinement with $245$
grid points at the final time (right). The accuracies are
comparable.}

\section{Summary}\label{sec:Summary}
We have developed an error control strategy for finite difference solutions
of parabolic equations, involving both temporal and spatial
discretization errors. The global time error strategy discussed in
\cite{LaVe2007} appears to provide an excellent starting point for
the development of such an algorithm. The classical ODE approach
used there and the principle of tolerance proportionality are combined with an efficient estimation
of the spatial error and mesh adaptation to control the overall
global error. Two approaches have been presented to handle spatial
mesh improvement: (i) globally uniform refinement and (ii) local refinement
and coarsening based on an equidistribution principle. Inspired by \cite{Be1988},
we have used Richardson extrapolation to approximate the spatial truncation
error within the method of lines. Our control strategy aims at balancing the spatial
and temporal discretization error in order to achieve an accuracy imposed
by the user.

The key ingredients are: (i) linearized error transport equations
equipped with sufficiently accurate defects to approximate the
global time error and global spatial error and (ii) uniform or
adaptive mesh
refinement and local error control in time based on tolerance
proportionality to achieve global error control. For illustration of
the performance and effectiveness of our approach, we have
implemented second-order finite differences in one space dimension and the
example integrator ROS3P \cite{LaVe2001}. On the basis of three
different test problems we could observe that our approach is very
reliable, both with respect to estimation and control.

Needless to say that spatial mesh adaptation locally in time is
more efficient for solutions having a strongly nonuniform nature
in space, especially if it varies over time. This is clearly visible
for the travelling wave solution of the Allen-Cahn problem.
However, optimized uniform strategies might also be of interest if
users would like to extend their own software packages not having
the option of dynamic adaptive mesh refinement to global error control.

%\newpage
%\input{Reviewer1.tex}
%\input{Review_CMAME-D-14-00245.tex}
%\input{Reviewer3.tex}
\end{document}